\documentclass[11pt]{article}
\usepackage{graphicx}
\newtheorem{theorem}{Theorem}
\newtheorem{lemma}[theorem]{Lemma}

\newtheorem{proposition}[theorem]{Proposition}

\newtheorem{definition}[theorem]{Definition}

\newtheorem{notation}[theorem]{Notation}

\begin{document}

\title{On dissimilarity vectors
of (not necessarily positive) weighted trees}
\author{  Elena Rubei }
\date{\hspace*{1cm}}
\maketitle

\vspace{-1.1cm}

{\small 
{\bf Address:}
Dipartimento di Matematica ``U. Dini'', 
viale Morgagni 67/A,
50134  Firenze, Italia 

{\bf
E-mail address:}
 rubei@math.unifi.it}

\bigskip

\smallskip

\def\thefootnote{}
\footnotetext{ \hspace*{-0.36cm}

}

{\small {\bf Abstract.} 
Let $T$ be  a (not necessarily positive) weighted tree  
 with $n$ leaves numbered by the set $\{1,...,n\}$.
Define the $k$-weights of the tree
 $D_{i_1,....,i_k}(T)$ as the sum of the lengths of the edges of the 
minimal subtree connecting $i_1$,....,$i_k$. 
We will call such numbers ``$k$-weights'' of the tree.
In this paper, 
 we characterize  sets of real numbers indexed by the  
subsets of any cardinality  
$ \geq 2$ of a $n$-set  to be the weights of a tree with $n$ leaves. }

\bigskip

\section{Introduction}

Consider a positive-weighted tree $T$ (that is a tree such that every 
edge is endowed with a positive real number, which we  call the 
length of the edge) with $n$ leaves numbered by the set 
$\{1,...,n\}$. 
Let $D_{i,j}(T)$ be the sum of the lengths of the edges of the shortest 
path connecting $i$ and $j$ for any $i$ and $j$ leaves of $T$. 
We call such number the ``double weight'' for $i$ and $j$. 

In 1971 Buneman characterized the metrics on finite sets which are the double
weights of a positive-weighted tree:

\begin{theorem} {\bf (Buneman)} A metric $(D_{i,j})$ 
on $\{1,...,n\}$ is the metric 
induced
 by a positive-weighted tree if and only if for all $i,j,k,h  \in \{1,...,n\}$
the maximum of $\{D_{i,j} + D_{k,h},D_{i,k} + D_{j,h},D_{i,h} + D_{k,j}
 \}$ 
is attained at least twice. 
\end{theorem}

The problem of reconstructing trees from data involving the distances  
between the leaves has several applications, such as phylogenetics:
 evolution of species can be represented by trees and, given
 distances between genetic  sequences of some species, one 
can  try to reconstruct the  evolution tree from these distances.
Some algorithms to reconstruct trees from the data $\{D_{i,j}\}$ have been 
proposed. Among them is neighbour-joining method, invented by
 Saitou and Nei in 1987 (see \cite{NS},
\cite{SK} and \cite{PSt2}).

For any weighted  tree $T$ with leaves $1,...,n$ and for any distinct 
$ i_1,..., i_k \in \{1,..., n\}$, 
define  $D_{i_1,....,i_k}(T)$ as the sum of the lengths of the edges 
of the minimal  subtree connecting $i_1$,....,$i_k$. We call such numbers 
$k$-weights of the tree $T$ and the vector of the $k$-weights is called 
$k$-dissimilarity vector.

In 2004, Pachter and Speyer proved the following theorem (see \cite{PS}).

\begin{theorem} {\bf (Pachter-Speyer)}. Let $ k ,n  \in {\bf N}$ with
$ n \geq 2k-1$ and $ k \geq 3$.  A positive-weighted tree
 $T$ with $n$ leaves $1,...,n$ and no vertices of degree 2
is determined by the values $D_I$ where $ I $ varies in the $k$-subsets 
of $\{1,...,n\}$.
\end{theorem}

It can be interesting to characterize the sets of real numbers 
which are sets of $k$-weights of a tree. 
% in fact (I quote Speyer and Sturmfels's  paper \cite{SS2})
%``it can be more reliable statistically
%to estimate the triple weights $D_{i,j,k}$ rather than the pairwise
%distances $D_{i,j}$''.  We refer to \cite{SS2} and above all to 
%\cite{PS} for an analysis of this and the references.

In \cite{Iri},  Iriarte proves that $k$-dissimilarity vectors  
 of positive-weighted  trees are contained in the tropical Grassmannian.
See also \cite{Co} and \cite{Man}.

Less results are known about not necessarily positive weighted trees.
Observe that also in this case, 
the problem of reconstructing the weighted trees may have some applications:
imagine that a particle, by going through an edge 
of a tree, gets or looses some substance (as much as the weight of the
 edge). 
If we know how much the substance of this particle varies by going from a 
leaf $i$ of the tree to another leaf $j$ (the value 
$D_{i,j}$) for any $i$ and $j$, we can try to reconstruct the weighted tree
(which can repesents a tree in the human body, a hydraulic web...).
Analogously the numbers $D_{i_1,...,i_k}$  
can represent how much a material,
 by going from $i_s$ to $i_1,...., \hat{i_s}, ..., i_k $, gets or 
looses of a certain substance.
It can be interesting, given a set 
$\{D_{i_1,...,i_k}\}_{i_1,...,i_k}$, to wonder if there exists 
a weighted tree with these $k$-weights.

In \cite{Ru} we gave a characterization for sets indexed by 
$2$-subsets (or $3$-subsets) of a $n$-set to be double (resp. triple) 
weights of a tree with $n$ leaves
 (with  not necessarily positive weights) 
and, by using these ideas, we proposed 
 a slight modification of Saitou-Nei's Neighbour-Joining 
algorithm to reconstruct trees from the data $D_{i,j} $.

Here we characterize sets of real numbers indexed by 
the subsets of any cardinality $ \geq 2$ 
of a $n$-set to be  the weights of a tree (Theorem
\ref{anyk}); besides
we extend the definition of $D_{i_1,....,i_k}(T)$ to the case 
$i_1,...., i_k$ not distinct and we find  necessary and sufficient 
conditions for a set of real numbers  indexed by  the submultisets of 
an $n$-set to be the set of the weights of a tree  with $n$ leaves 
(Theorem \ref{anykmulti}) and   necessary and sufficient 
conditions for a set of real numbers  indexed by  the $k$-submultisets of 
an $n$-set to be the set of the $k$-weights of a tree  with $n$ leaves 
(Theorem \ref{kweights}).

%NOTATION

\section{Some notation}

\begin{definition} \label{bell}
A {\bf cherry}  $B$ in a tree $T$ is  a subtree  such that only one of 
the inner vertices is not bivalent;
we call this vertex  ``{\bf stalk}'' of the cherry
and we say that the  leaves of $B$  are {\bf neighbours}. 
We call ``{\bf twig}'' of a leaf of a cherry 
the path from the leaf to the stalk of the cherry.

A {\bf complete cherry} is a cherry such that there doesn't exist another 
cherry strictly containg it.

\end{definition}

\begin{center}
\includegraphics[scale=0.45]{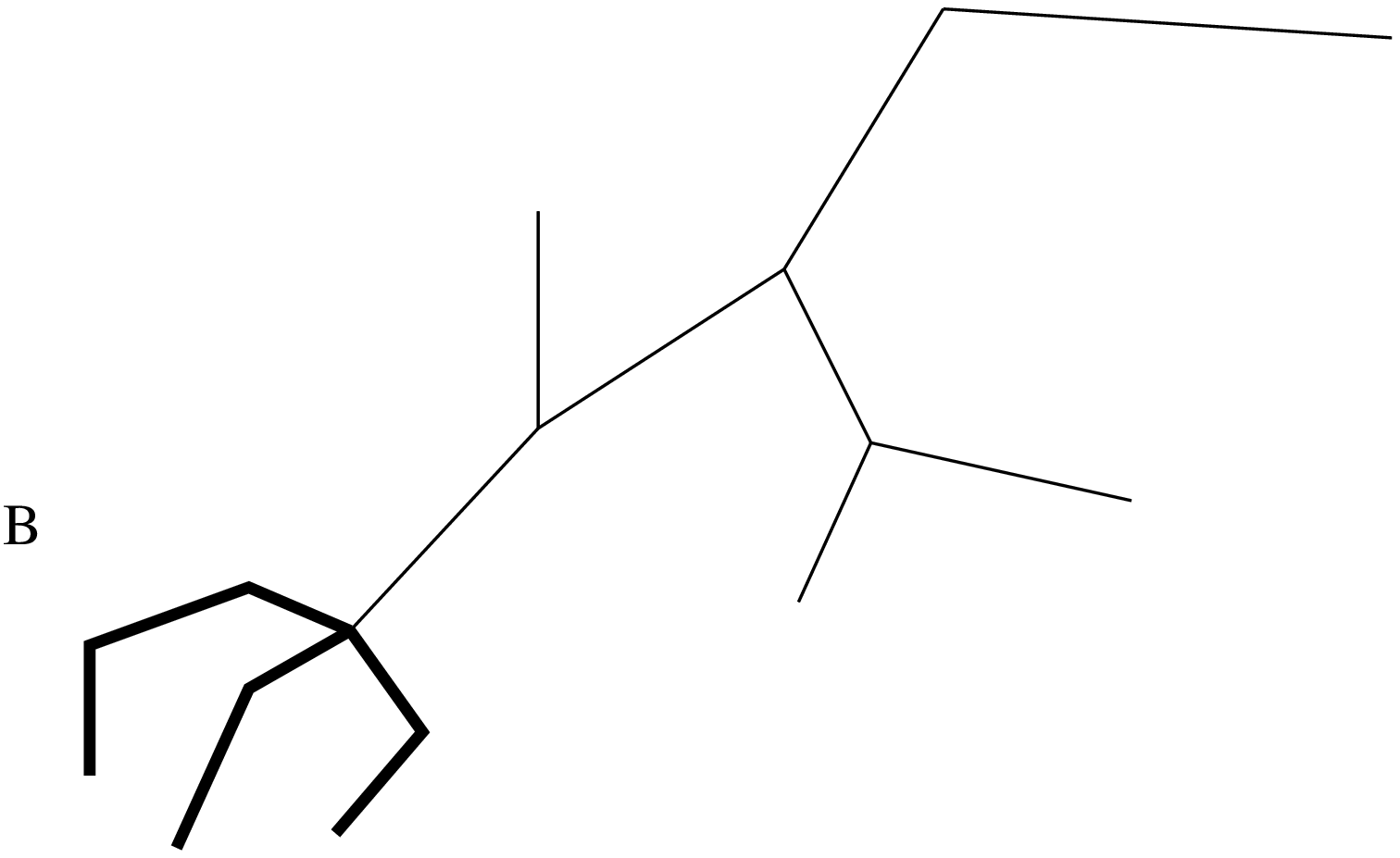}
\end{center}

\begin{notation} \label{trees}
$\bullet$ For every $n \in {\bf N}$, let $ [n] = \{1,....,n\}$. 

 $\bullet$ If $M$ is a set,  
 denote  the set of the $k$-subsets (without repetitions) 
of $M$ by $M^k$ and  
the set of the ``$k$-submultisets'' (i.e. with possible repetitions)
of $M$ by $M_k$. 

\smallskip

$\bullet$ 
 For any  set $\{ D_{\{i_1,...., i_k\}}\}$ of real numbers indexed by 
 elements  $ \{i_1,...., i_k\} \in [n]^k$ or $[n]_k$,
we  denote $D_{\{i_1,..., i_k\}}$ by $ D_{i_1,..., i_k}$ 
for any order of $i_1,...., i_k$.

\smallskip

$\bullet$ A {\bf weighted tree} is a tree such that every edge is 
endowed with a real number called weight or length of the edge. If 
the weights are positive we say that the tree is {\bf positive-weighted}.
Please note that in other papers ``weighted'' means positive-weighted.

For $x, y$ vertices of a tree $T$, we denote by $d(x,y) 
$ ({\bf intrinsic distance}) 
the  number of the edges of the path from $x$ to $y$.  

For any leaf $x$ of $T$ and any subtree $E$, we define 
$ {\bf N(x,E)}$ as the at least trivalent vertex in $E$
 with minimum intrinsic distance from $x$. 

\medskip

{\it Example.}
\begin{center}
\includegraphics[scale=0.42]{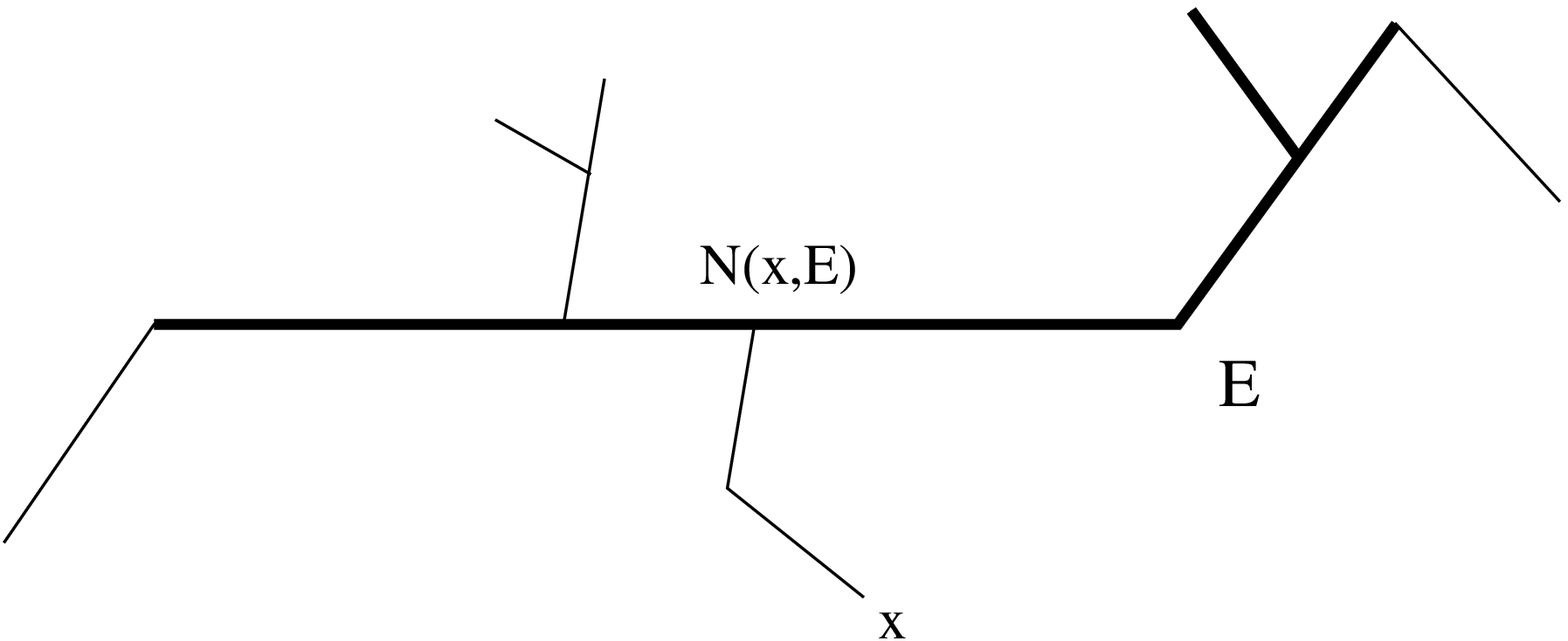}
\end{center}

%For any subset $S$ of leaves of $T$, we denote by $\langle S \rangle$ the
% minimum subtree of $T$ containg $S$.

Now let $T$ be weighted and let $[n]$ be the set of its leaves.
For $x, y$ vertices of $T$, we denote by $w(x,y) $ ({\bf w-distance}) 
the sum of the weights 
of the edges of the path from $x$ to $y$ (obviously it is not a distance).

For any distinct $i_1,..., i_k \in [n]$, we define $D_{i_1, ..., i_k} (T)$
as the sum of the lengths of the edges 
of the minimal  subtree connecting $i_1$,....,$i_k$. Besides if $ 
i_1,...,  i_k $ are not distinct, we define  $D_{i_1, ..., i_k} (T)$ 
by induction on the number of the repetition in the following way:
$$ D_{x, x, Z} (T) = D_{x,Z} (T) +  w(x, N(x,T))$$
We call the numbers $D_{i_1,..., i_k}(T)$ {\bf $k$-weights of $T$}.
\end{notation}

{\it Example.}
\begin{center}
\includegraphics[scale=0.45]{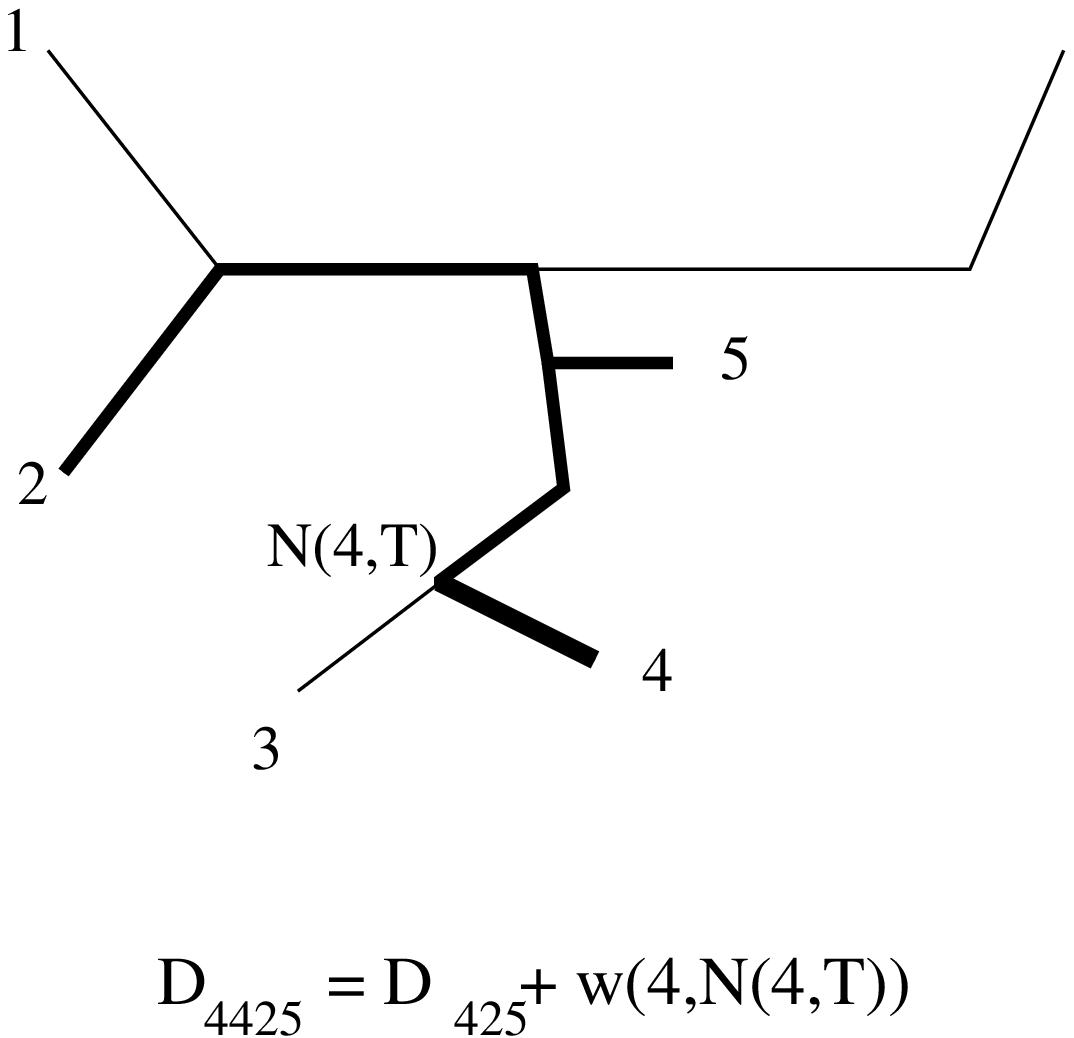}
\end{center}

\begin{definition} \label{pseudobell}
For any  set $\{ D_I \}_{I \in  S}$ of real numbers parametrized
by $S \subset {\cal M}:=\{I \;submultisubset\; of [n]\}$  and for any 
$e, e' \in [n]$, we  define $\ast^{e, e'} $ the 
following condition:
\begin{center}
$  D_{e, X}  - D_{e', X } $
doesn't depend on $X $ s.t. $(e,X), (e',X) \in S$
\end{center}

and
we say that $\alpha = \{\alpha_1,....,\alpha_r\} \subset  
 [n]$ is a {\bf pseudocherry}  for $\{D_I\}_{I \in S}$ 
if $ \ast^{\alpha_i, \alpha_j}$ holds
for all $ i,j $. 
We say that $\alpha$ is a   {\bf complete pseudocherry}
if $\not \exists \gamma \in [n]- \alpha$
 such that $ \ast^{\gamma, \alpha_i}$  holds for all  $ i $. 
\end{definition}

\begin{proposition} \label{starbell} \cite{Ru}.
Let $T$ be a positive-weighted tree with leaves  $1,....,n$ with 
$n \geq 2k-1$.  Let $ e , e' \in [n]$.
Then $  \ast^{e, e'} $ holds for $ \{ D_I (T) \}_{I \in [n]^k}$ 
if and only if $ \{e ,e'\} $ is a cherry, that 
is  $ \{e ,e'\} $ is a a pseudocherry if and only if 
it is a cherry.
\end{proposition}

% TREES WITH FOUR LEAVES
\section{The case of the trees with four leaves}

\begin{lemma} \label{caso4}  Let $ k \in {\bf N}$, $ k
\geq 2$. 
Let $\{ D_{I}\}_{\{I \} \in [4]_k}$ be a set of  real numbers. 
It is the set of $k$-weights of a tree $T$ with $1,2,3,4$ as leaves and
$ \{1,2\}$ and $ \{3,4\}$ as 
cherries if and only if   the following conditions hold:

A)  $\{1,2\}$ and $\{3,4\} $ are pseudocherries and
 for any $ i \in \{1,2\} $, $ j \in \{3,4\} $
$$ D_{i , X} - D_{j, X} $$ 
is the same  for $X$ varying in the subsets of $[4]$ intersecting both
$\{1,2\} $ and $ \{3,4\}$, it is the same for $X$ varying in the subsets of 
$\{1,2\}$,  it is the same for $X$ varying in the subsets of 
$\{3,4\}$.

%$ \ast^{\delta, \delta'}_{<} $ holds 
%for any $ \delta, \delta' \in [4]$ for the order $ \{1,2\} < \{3,4\}$ and
%$ 1 \mid 2 $, $ 3 \mid 4$ and 

B) $$  -D_{i,Z} + D_{j,Z} - D_{i,Y} + D_{j,Y} = 2 (D_{j,W} - D_{i,W}) $$
for $ i \in \{1,2\} $ and $ j \in \{3,4\}$, 
for $Z \subset \{1,2\}$, $ Y \subset \{3,4\}$, 
$W$ intersecting both $ \{1,2\} $ and $ \{3,4\}$.  
\end{lemma}

{\it Proof.} $ \Rightarrow$ Easy.

$ \Leftarrow$ 
We will sometimes denote $D_{1,...,1,2,...,2,3,...,3,4,...,4}$, with 
$1$ repeated $k_1$ times, 
$2$ repeated $k_2$ times, 
$3$ repeated $k_3$ times, 
$4$ repeated $k_4$ times, 
 by $ D_{1^{k_1}, 2^{k_2}, 3^{k_3}, 4^{k_4}}$. 

%Observe that $  \ast^{\delta, \delta'}_{<} $ implies $ \ast^{1,2}$ and 
%$ \ast^{3,4}$.

\medskip

{\it Remark.} {\em 
If there exists a tree $T$ without bivalent vertices, with $[4]$ as set of 
leaves  and  $ \{1,2\}$ and $ \{3,4\}$  as cherries, 
 by calling the weights of the edges as in the figure below, we have:
 
if $i$ and $j$ are in the same cherry, then
$ a_i -a_j = D_{i,X} -D_{j,X}$ for any $ X \in [4]_{k-1}$ 

if $i$ and $j$ are not in the same cherry, then:

$ a_i -a_j = D_{i,X} -D_{j,X}$ for any $  X \in [4]_{k-1}$ 
intersecting both cherries and 

$ a_i +f -a_j = D_{i,X} -D_{j,X}$ for any $  X \in [4]_{k-1}$ in
the same cherry as $j$.}

\begin{center}
\includegraphics[scale=0.4]{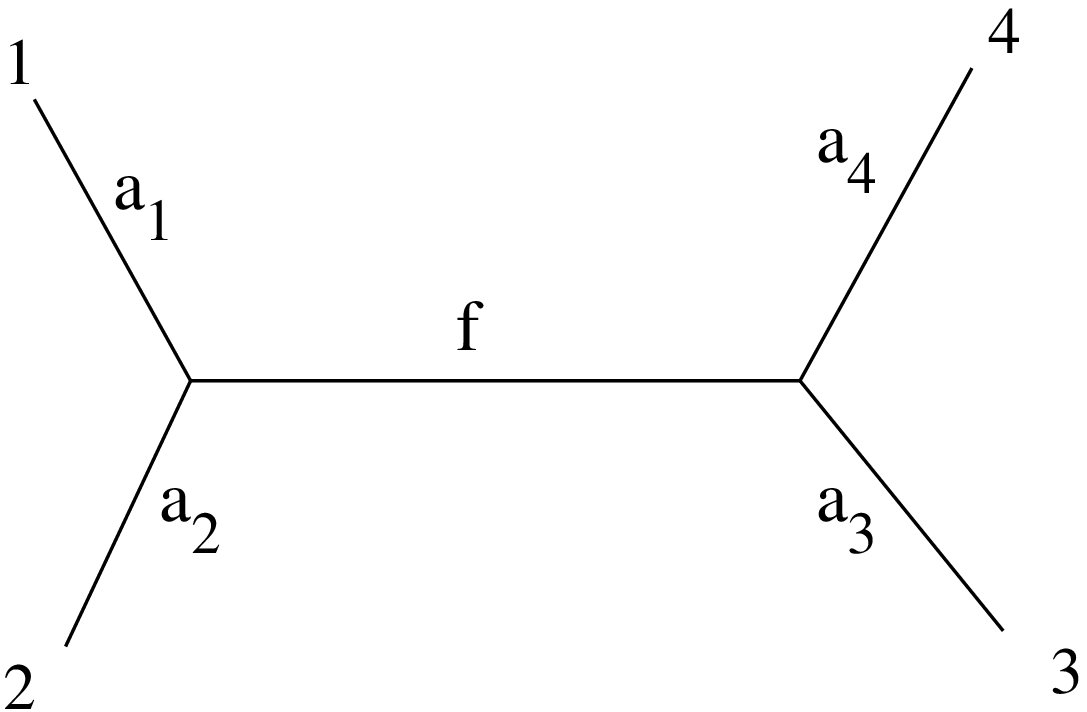}
\end{center}

We define a tree as in the figure above with $ a_1,a_2,a_3,a_4,f$ defined 
in the following way.

We define $a_1, a_2 $ as the  solution of the following linear system 
(for any $k_1, k_2$ with $ k_1 + k_2=k$ and any $X \in [4]_{k-1}$):

$ \left\{ \begin{array}{l} a_1 -a_2 = D_{1,X} -D_{2,X} \\
k_1 a_1 + k_2 a_2 = D_{1^{k_1}, 2^{k_2}, 3^0, 4^0}
\end{array}\right.$ 

Obviously it admits only one solution and one can easily see that 
it does not depend on $X$ since $ \ast^{1,2}$ holds; besides 
it doesn't depend  on $ k_1, k_2$: in fact the system

$ \left\{ \begin{array}{l} a_1 -a_2 = D_{1,X} -D_{2,X} \\
k_1 a_1 + k_2 a_2 =  D_{1^{k_1}, 2^{k_2}, 3^0, 4^0} \\
t_1 a_1 + t_2 a_2 =  D_{1^{t_1}, 2^{t_2}, 3^0, 4^0}
\end{array}\right.$ 

is compatible for any $t_1, t_2$ with
$ t_1 +t_2 =k$; to see this, it is sufficient to see 
that 

$ \left\{ \begin{array}{l} a_1 -a_2 = D_{1,X} -D_{2,X} \\
k_1 a_1 + k_2 a_2 =  D_{1^{k_1}, 2^{k_2}, 3^0, 4^0} \\
(k_1-1) a_1 + (k_2+1) a_2 =  D_{1^{k_1 -1}, 2^{k_2 +1}, 3^0, 4^0}
\end{array}\right.$ 

is compatible and this follows from $ \ast^{1,2}$. 

We define: 
$$f= D_{i,Y} -D_{j,Y} -D_{i,W} +D_{j,W}$$
for any $ i,j$ with $i 
\in \{1,2\} $, $ j \in \{3,4\}$ and for $W $ intersecting both 
$\{1,2\}$ and $ \{3,4\}$ and $ Y $ in $ \{3,4\}$.
One can easily see that it is a good definition, that is  it does not
 depend on $ i,j, W, Y$. 
Equivalently (by $B$)  we can define
$$ f= -D_{i,Z} +D_{j,Z} +D_{i,W} -D_{j,W}$$
for any $ i,j$ with $i 
\in \{1,2\} $, $ j \in \{3,4\}$ and for $W $ intersecting both 
$\{1,2\}$ and $ \{3,4\}$ and $ Z $ in $ \{1,2\}$.

For  $ j \in \{3,4\}$ we define 
$$ a_j =  a_1 + (D_{j,W} -D_{1,W}) $$
for $W$  intersecting both 
$\{1,2\}$ and $ \{3,4\}$. It is a good definition by $A$. 

%\smallskip

%{\it Remark.} {\em With the above definitions, we have that 
%$$ a_i -a_j = D_{i,X}-D_{j,X}$$ 
%for any $X$ if $ i$ and $j$ are both in $ \{1,2\}$ or in $\{3,4\}$, 
%for $X$ intersecting both $ \{1,2\}$ and $ \{3,4\}$  if $i$ is one 
%of the two set and $j$ in the other.}

\smallskip

We shall show now that  for such a tree $T$ we have 
$$   D_{1^{k_1}, 2^{k_2}, 3^{k_3}, 4^{k_4}} (T)=
 D_{1^{k_1}, 2^{k_2}, 3^{k_3}, 4^{k_4}}$$ 
$ \bullet $ Let us first suppose 
 that $ k_1 + k_2 > 0$. We argue by induction
on $ k_3 + k_4 $. 

If $ k_3 + k_4 =0 $ it is obvious by the definition of $ a_1$ and $ a_2$:
$$   D_{1^{k_1}, 2^{k_2}, 3^{0}, 4^{0}} (T)  
\stackrel{de\!f \; o\!f \; T }{=} k_1 a_1 + k_2 a_2 
 \stackrel{  de\!f \; o\!f \;a_1 \; and \; a_2 }{=}  
D_{1^{k_1}, 2^{k_2}, 3^{0}, 4^{0}} $$
If $ k_3 + k_4 =1 $, we can suppose for instance that $k_3 =1$ and $ k_4=0 $.
$$   D_{1^{k_1}, 2^{k_2}, 3^{1}, 4^{0}} (T)=
 D_{1^{k_1+1}, 2^{k_2}, 3^{0}, 4^{0}} (T) +f + a_3 -a_1=$$
$$ \stackrel{by \;induct.\; assumpt. \;and\;2^{nd}\; 
de\!f. \; o\!f \; f \; and \; a_3}{=}
D_{1^{k_1+1}, 2^{k_2}, 3^{0}, 4^{0}}  +D_{1^{k_1}, 2^{k_2}, 3^{1}, 4^{0}} -
 D_{1^{k_1+1}, 2^{k_2}, 3^{0}, 4^{0}} 
= D_{1^{k_1}, 2^{k_2}, 3^{1}, 4^{0}} $$ 
If $ k_3 + k_4 -1>0 $, 
$$   D_{1^{k_1}, 2^{k_2}, 3^{k_3}, 4^{k_4}} (T)=
 D_{1^{k_1+1}, 2^{k_2}, 3^{k_3-1}, 4^{k_4}} (T)  + a_3 -a_1=$$
$$\stackrel{by \; induct.\; assumpt. \; and\;de\!f. \; o\!f \; a_3}{=}
D_{1^{k_1+1}, 2^{k_2}, 3^{k_3-1}, 4^{k_4}} +
  D_{1^{k_1}, 2^{k_2}, 3^{k_3}, 4^{k_4}}- 
D_{1^{k_1+1}, 2^{k_2}, 3^{k_3-1}, 4^{k_4}} = 
 D_{1^{k_1}, 2^{k_2}, 3^{k_3}, 4^{k_4}}
 $$ 
$\bullet $ Suppose now that $ k_1 + k_2 =0$.  
$$ D_{1^{0}, 2^{0}, 3^{k_3}, 4^{k_4}} (T) 
=  D_{1^1, 2^{0}, 3^{k_3-1}, 4^{k_4}}(T) -f -a_1 + a_3 = $$ 
$$ \stackrel{by \;previous\; case 
\;and\;1^{st} \; de\!f. \; o\!f \; f \; and \; a_3 
}{=} 
  D_{1^1, 2^0, 3^{k_3-1}, 4^{k_4}} + 
  D_{1^0, 2^0 , 3^{k_3}, 4^{k_4}}
- D_{1^1, 2^0 , 3^{k_3-1}, 4^{k_4}}{=}
  D_{1^0, 2^0, 3^{k_3}, 4^{k_4}}
$$

\hfill \framebox(7,7)

\begin{lemma} \label{caso4senzarip}
Let $\{ D_{I}\}_{\{I  \subset [4] \;|\; cardinality (I) \geq 2\}}$ 
be a set of  real numbers. 
It is the set of weights of a tree $T$ with $1,2,3,4$ as leaves and
$ \{1,2\}$ and $ \{3,4\}$ as 
cherries if and only if  $\{1,2\}$ and $ \{3,4\}$ are pseudocherries and 

A) $$  -D_{i,Z} + D_{j,Z} - D_{i,Y} + D_{j,Y} = 2 (D_{j,W} - D_{i,W}) $$
for $ i \in \{1,2\} $ and $ j \in \{3,4\}$, 
for $Z \subset \{1,2\}$, $ Y \subset \{3,4\}$, 
$W$ intersecting both $ \{1,2\} $ and $ \{3,4\}$.  

B) if we define $a_{i}$ for any $ i \in \{1,2\} $ by 
$ a_{i}  = 
\frac{1}{2}( D_{i, j} + D_{i,X} -D_{j,X}) $ 
 for any $X \subset [n]$ and $ j \in \{1,2 \} -\{i\}$, then,
for $ \{i,j \} = \{1,2\} $ and any $ \delta \subset [n]$, 
$$D_{i,j, \delta} = 
  a_i +  D_{j, \delta} $$

\end{lemma}

{\em Sketch of the proof.}
$\Rightarrow$ Easy to prove.

$\Leftarrow$  Construct a tree $T$ as in the proof of Lemma \ref{caso4}
from the $ D_{I}$ with $cardinality(I) =2$. We have to prove that 
$ D_{I}(T) = D_I $ for any $ I \subset [n]$ with $cardinality (I) =3,4$.
$$ D_{1,2,3} (T) = a_1 + D_{2,3} (T) =   a_1 + D_{2,3} =  D_{1,2,3}  $$ 
Analougously for $ D_{1,2,4} $. 
$$  D_{2,3,4}(T) = D_{1,2,4} (T)+\frac{1}{2}(
 -D_{1,2} (T) + D_{3,2} (T) - D_{1,4}(T)
+ D_{3,4} (T)) =  D_{1,2,4} +\frac{1}{2}( -D_{1,2}  + D_{3,2}  - D_{1,4}
+ D_{3,4} ) \stackrel{A}{=}  D_{2,3,4}$$
Analogously $ D_{1,3,4}$. 
$$ D_{1,2,3,4} (T) = a_1 + D_{2,3,4} (T) =   a_1 + D_{2,3,4} 
=  D_{1,2,3,4}  $$ 
 
\hfill \framebox(7,7)

% CHARCATERIZATION...
\section{Characterization of the set of dissimilarity  vectors}

In this section our first aim is to characterize 
the sets of real numbers indexed by subsets or submultisets of $ [n]$ 
which come from a tree.  We  characterize also
the sets of real numbers indexed by
the elements of  $[n]_k$ for $k$ fixed. 
Shortly speaking in \cite{Ru} we proved that for $k=2$ 
such a set comes from a tree if and only if in 
$[n]$ there are at least two pseudocherries and if we substitute
every pseudocherry with a point, the same condition holds for the new set
 and so on. Obviously for higher $k$ the situation is more complicated. 
%because, while for $k=2$ we can recover the weight of a twig of a cherry 
%$i,j$ only from $D_{i,j}$ and $ D_{i,x}- D_{j,x}$ for any $x$, in case 
%of higher $k$, to recover the length of a twig  we need more data. 

\begin{theorem}  \label{anykmulti} Let $n \geq 4$. 
Let $\{ D_{I}\}_{\{I  \in [n]_k \; for \; some \; k \geq 2\}}$ 
be a set of real numbers. 
It is the set of the  weights of a tree $T$ with leaves $ 1,...,n$ 
if and only if there exist  $\alpha, \beta  \subset [n]$ such that:

1) $\alpha $ and $ \beta $ are  disjoint complete pseudocherries and
for any $ \alpha_i \in \alpha$, $ \beta_j \in \beta $ the number
$$ D_{\alpha_i , X} - D_{\beta_j, X}$$ 
is the same  for $X$ varying in the submultisets of $[n]$ intersecting both
$\alpha $ and $ \beta$, it is the same for $X$ varying in the submultisets of 
$\alpha$,  it is the same for $X$ varying in the submultisets of 
$\beta$

2) $$ -D_{\alpha_i,Z} + D_{\beta_j,Z} - D_{\alpha_i,Y} + D_{\beta_j,Y} = 
2 (D_{\beta_j,W} - D_{\alpha_i,W}) $$
for $\alpha_i \in \alpha $, $ \beta_j \in \beta$, 
$Z \subset \alpha$, $ Y \subset \beta $ and 
$W$ intersecting both $ \alpha $ and $ \beta$ 

3) if, for any $ \alpha_i \in \alpha $, we define $a_{\alpha_i}$  by 
$ a_{\alpha_i}  = 
\frac{1}{2}( D_{\alpha_i, \alpha_j} + D_{\alpha_i,X} -D_{\alpha_j,X}) $ 
 for any $X$ submultiset of $ [n]$ and $ \alpha_j \in \alpha$, then,
for any $ \alpha_1,..., \alpha_{t} 
\in \alpha  $ and $ \delta_1,..., \delta_s \in [n]$, 
$$D_{\alpha_1,..., \alpha_t, \delta_1,..., \delta_s} = 
  a_{\alpha_1}+...+ a_{\alpha_{t-1}} +
  D_{\alpha_t, \delta_1,..., \delta_s} $$  
(i.e.
$2( D_{\alpha_1,..., \alpha_t, \delta_1,..., \delta_s}-
 D_{\alpha_t, \delta_1,..., \delta_s}) = 
  D_{\alpha_1, \alpha_2}+ D_{\alpha_2, \alpha_3}+...+ 
D_{\alpha_{t-2}, \alpha_{t-1}} + D_{\alpha_1, X} -D_{\alpha_{t-1},X} $  
 $\forall X \subset [n]$)

4)  if we define  $ M = [n] - \alpha \cup \{\underline{\alpha}\} $
and $$ D_{\underline{\alpha} , i_1,..., i_{k}} = 
D_{\alpha_i , i_1,..., i_{k}} -a_{\alpha_i}$$ 
(for any $ \alpha_i \in \alpha$), 
then the same conditions  hold for $M$.
 
\end{theorem}

{\em Proof.}
$\Rightarrow$ Easy to prove; for instance 
observe that, for $X$ subset of $ \alpha$, 
 $D_{\alpha_i,X} - D_{\beta_j, X} $
is the sum of the weights of the edges of the twig of $\alpha_i$, 
minus the sum of the  weights of the edges of the twig of $\beta_j$, 
minus the w-distance between the stalks of $\alpha $ and $ \beta$.

\begin{center}
\includegraphics[scale=0.4]{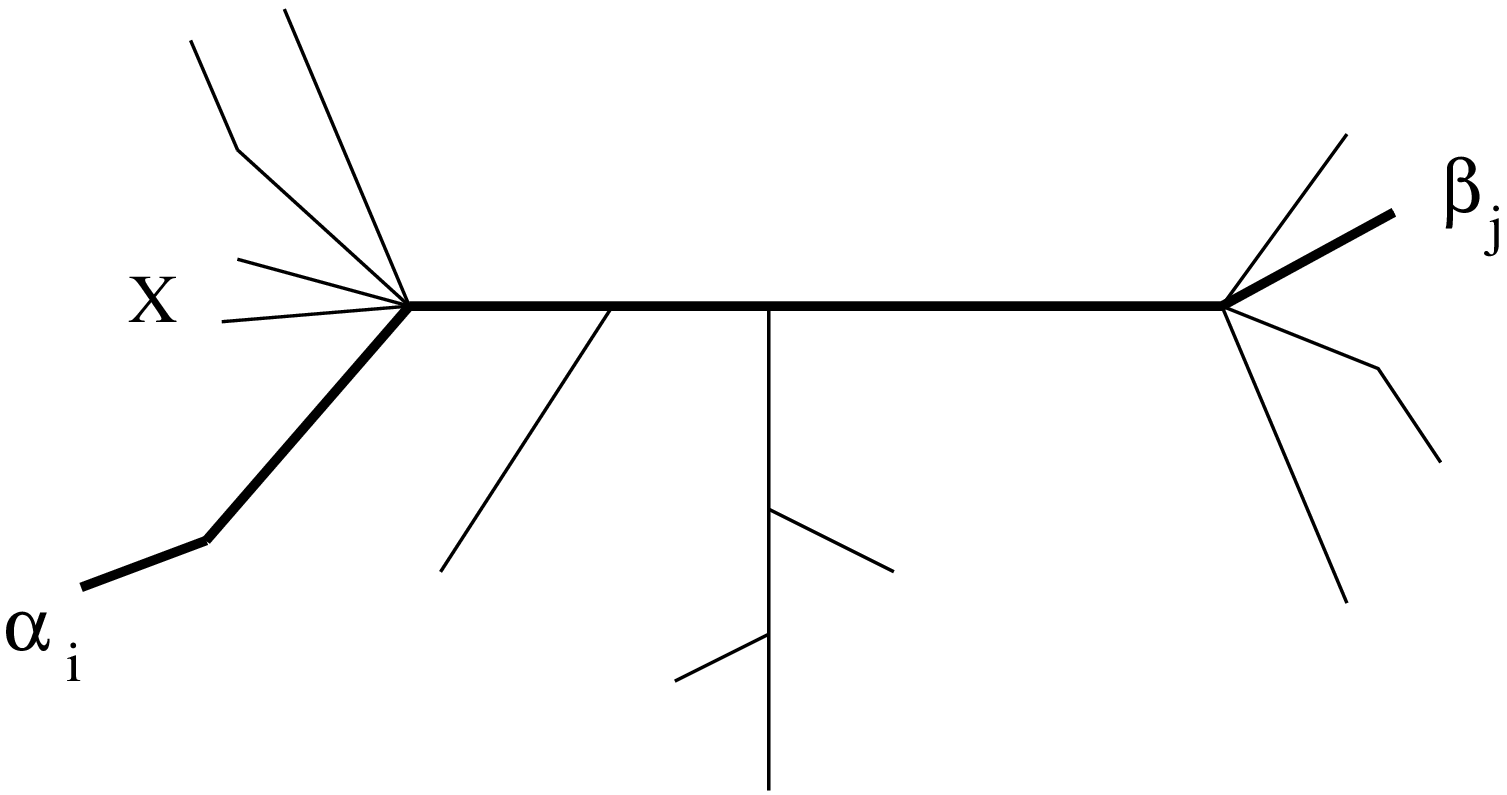}
\end{center}

$\Leftarrow$ 
First observe  that the definition of $ a_{\alpha_i}$ 
is equivalent to the formula

$ \left\{ \begin{array}{l} a_{\alpha_i} -a_{\alpha_j} = D_{\alpha_i,X} -
D_{\alpha_j,X} \\
 a_{\alpha_i} + a_{\alpha_j} = D_{\alpha_i, \alpha_j} 
\end{array}\right.$
  
 for any $X$  submultiset of $[n]$ and it
is a  good definition, that is it doesn't depend on $X$ neither on 
$ \alpha_j$  (because $ \alpha$ is a pseudocherry). Besides, obviously, 
also the  definition of $  D_{\underline{\alpha} , i_1,..., i_{k}} $ 
doesn't depend on  $\alpha_i$, because $ \alpha$ is a pseudocherry.

We can prove the statement by induction on $n$.
The case $n=4$ follows from Lemma \ref{caso4} 
(observe that conditions 1 and 2 of the theorem imply conditions A and B 
of the lemma and that the definitions of $a_i$ and $f$ in the proof 
of the lemma  don't depend on $k$, so the tree we construct is the 
same for every $k$). 

Let us prove the induction step.
By induction assumption there exists  a tree $R$ such that  
$$ D_{I}(R)= D_{I } \;\;\; \forall I  \in M_k \; for \; k \geq 2$$ 
We define the tree $T$ by attaching to $ R$  a cherry $ \alpha$ with 
lengths of the twigs $ a_{\alpha_i}$ to the point $\underline{\alpha}$.
 We must show that 
$ D_{\alpha_1,..., \alpha_t, \delta_1,...., \delta_{s}} (T) = 
D_{\alpha_1,..., \alpha_t, \delta_1,...., \delta_{s}} $
for any $\alpha_1, ..., \alpha_t \in \alpha$,
 $ \delta_1,...., \delta_{s}  \in [n] -\alpha $. 
We prove this on induction on $t$. 

$t=0$ is obvious 

\medskip

$ t=1$  $$ D_{\alpha_j,\delta_1,...., \delta_{s}} (T) 
\stackrel{de\!f.\; o\!f \; T}{=} 
a_{\alpha_j} +D_{\underline{\alpha}, \delta_1,...., \delta_{s}} (R) 
\stackrel{ind. \; ass.}{=}
a_{\alpha_j} +D_{\underline{\alpha}, \delta_1,...., \delta_{s}} 
\stackrel{de\!f.\; of \;D_{\underline{\alpha} , \delta_1,...., 
\delta_{s}}}{=}$$ 
$$ = a_{\alpha_j} +D_{\alpha_j, \delta_1,...., \delta_{s}} -a_{\alpha_j} 
=D_{\alpha_j, \delta_1,...., \delta_{s}} $$ 

\medskip

Induction step:
$$ D_{\alpha_1,..., \alpha_t,\delta_1,...., \delta_{s}} (T) 
\stackrel{de\!f.\; o\!f \; T}{=} 
a_{\alpha_1} +....+ a_{\alpha_t} +
D_{\underline{\alpha}, \delta_1,...., \delta_{s}} (R) 
\stackrel{ind. \; ass.}{=}$$ 
$$ = a_{\alpha_1} +....+ a_{\alpha_t} +
D_{\underline{\alpha}, \delta_1,...., \delta_{s}} 
\stackrel{de\!f.\; of \;D_{\underline{\alpha} , \delta_1,...., 
\delta_{s}}}{=} a_{\alpha_1} +....+ a_{\alpha_t} +
D_{\alpha_1, \delta_1,...., \delta_{s}} -a_{\alpha_1}= $$
$$ = a_{\alpha_2} +....+ a_{\alpha_t} +
D_{\alpha_1, \delta_1,...., \delta_{s}}
\stackrel{3}{=} D_{\alpha_1,..., \alpha_t,\delta_1,...., \delta_{s}} $$
\hfill \framebox(7,7)

\bigskip

With a completely analogous  proof, we can prove the following statement for 
a set of real numbers indexed by subsets of $[n]$ (without repetitions):

\bigskip

\begin{theorem}  \label{anyk} Let $n \geq 4$. 
Let $\{ D_{I}\}_{\{I  \subset [n] \;|\; cardinality(I) \geq 2\}}$ 
be a set of real numbers. 
 It is the set of the  weights of a tree $T$ with leaves $ 1,...,n$ 
if and only if there exist two disjoint 
complete pseudocherry $\alpha $ and $ \beta$ in $  [n]$ 
such that

1) $$ -D_{\alpha_i,Z} + D_{\beta_j,Z} - D_{\alpha_i,Y} + D_{\beta_j,Y} = 
2 (D_{\beta_j,W} - D_{\alpha_i,W}) $$
for $\alpha_i \in \alpha $, $ \beta_j \in \beta$, 
$Z \subset \alpha$, $ Y \subset \beta $ and 
$W$ intersecting both $ \alpha $ and $ \beta$.  

2) if,  for any $ \alpha_i \in \alpha $,
 we define $a_{\alpha_i}$ by 
$ a_{\alpha_i}  = 
\frac{1}{2}( D_{\alpha_i, \alpha_j} + D_{\alpha_i,X} -D_{\alpha_j,X}) $ 
 for any $X \subset [n]$ and any $ \alpha_j \in \alpha$, then,
for any $ \alpha_1,..., \alpha_{t} 
\in \alpha  $ and  $ \delta_1,..., \delta_s \in [n]$, 
$$D_{\alpha_1,..., \alpha_t, \delta_1,..., \delta_s} = 
  a_{\alpha_1}+...+ a_{\alpha_{t-1}} +
  D_{\alpha_t, \delta_1,..., \delta_s} $$  
(i.e.
$2( D_{\alpha_1,..., \alpha_t, \delta_1,..., \delta_s}-
 D_{\alpha_t, \delta_1,..., \delta_s}) = 
  D_{\alpha_1, \alpha_2}+ D_{\alpha_2, \alpha_3}+...+ 
D_{\alpha_{t-2}, \alpha_{t-1}} + D_{\alpha_1, X} -D_{\alpha_{t-1},X} $  
 $\forall X \subset [n]$)

3)  if we define  $ M = [n] - \alpha \cup \{\underline{\alpha}\} $
and $$ D_{\underline{\alpha} , i_1,..., i_{k}} = 
D_{\alpha_i , i_1,..., i_{k}} -a_{\alpha_i}$$ 
(for any $ \alpha_i \in \alpha$), 
then the same conditions  hold for $M$.
 
\end{theorem} 

Finally we consider the case of set of real numbers indexed 
by $k$-submultiset of $[n]$ ($k$ fixed).

\begin{theorem}  \label{kweights} Let $n \geq 4$ and $ k \in {\bf N}$. 
Let $\{ D_{I}\}_{\{I \} \in [n]_k}$ be a set of  real numbers. 
It is the set of the $k$-weights of a tree $T$ with leaves $ 1,..., n$ 
if and only if   there exist  
$\alpha, \beta  \subset [n]$ such that:

1) $\alpha $ and $ \beta $ are disjoint complete pseudocherries and
for any $ \alpha_i \in \alpha$, $ \beta_j \in \beta $ the number
$$ D_{\alpha_i , X} - D_{\beta_j, X}$$ 
is the same  for $X$ varying in the subsets of $[n]$ intersecting both
$\alpha $ and $ \beta$, it is the same for $X$ varying in the subsets of 
$\alpha$,  it is the same for $X$ varying in the subsets of 
$\beta$.

2) $$ -D_{\alpha_i,Z} + D_{\beta_j,Z} - D_{\alpha_i,Y} + D_{\beta_j,Y} = 
2 (D_{\beta_j,W} - D_{\alpha_i,W}) $$

for $\alpha_i \in \alpha $, $ \beta_j \in \beta$, 
$Z \subset \alpha$, $ Y \subset \beta $ and 
$W$ intersecting both $ \alpha $ and $ \beta$.  

3)  if we define $a_{\alpha_i}$ for any $ \alpha_i \in \alpha $ by 
$ a_{\alpha_i}  = \frac{1}{k}( k_j D_{{\alpha_i}^{k_i}, {\alpha_j}^{k_j}} 
+ D_{\alpha_i,X} -D_{\alpha_j,X}) $ 
  for any $X \subset [n]$, any $ \alpha_j \in \alpha$,  
any $k_i , k_j \in {\bf N}$ 
with $k_i + k_j=k$, and analogously $a_{\beta_j}$, then 
$$  a_{\alpha_i} - a_{\beta_j} = D_{\alpha_i, A,D} + D_{\beta_j, B,D} 
- D_{\delta , A,D} - D_{\delta , B,D} $$ 
for any $ A \subset \alpha$, $ B \subset \beta$, $ D \ni \delta$.

4) If we define $ M = [n] - \alpha \cup \{\underline{\alpha}\} $
and $$ D_{\underline{\alpha} , i_1,..., i_{k-1}} = 
D_{\alpha_i , i_1,..., i_{k-1}} -a_{\alpha_i}$$ 
(for any $ \alpha_i \in \alpha$), 
then the same conditions  hold for $M$.

($X,Y,Z,W,A,B,D$ of size such that the size of all the indices is $k$.)

\end{theorem}

{\it Proof.}
$\Rightarrow$ Easy to prove.

\medskip

$\Leftarrow$ As in the proof of Theorem  \ref{anyk},
the definition of $ a_{\alpha_i}$ 
(which is equivalent to the formula

$ \left\{ \begin{array}{l} a_{\alpha_i} -a_{\alpha_j} = D_{\alpha_i,X} -
D_{\alpha_j,X} \\
k_i a_{\alpha_i} + k_j a_{\alpha_j} = D_{\alpha_i^{k_i}, \alpha_j^{k_j}} 
\end{array}\right.$
  
 for any $X \subset [n]$ and  any $k_i , k_j \in {\bf N}$ 
with $k_i + k_j=k$)
 and  the  definition of $  D_{\underline{\alpha} , i_1,..., i_{k}} $ are good definitions.
%doesn't depend on  $\alpha_i$, because $ \alpha$ is a pseudocherry.

We can prove the statement by induction on $n$.
The case $n=4$ follows from Lemma \ref{caso4}. 
%(observe that conditions 1 and 2 of Theorem imply conditions A and B of Lemma).
Let us prove the induction step.
By induction assumption and condition 4, there exists  a tree $R$ such that  
$$ D_{I}(R)= D_{I } \;\;\; \forall I  \in M_k $$ 
We define the tree $T$ by attaching to $ R$  a cherry $ \alpha$ with 
lengths of the twigs $ a_{\alpha_i}$ to the point $\underline{\alpha}$.
 We must show that for any $ \alpha_1, ..., \alpha_{t} \in \alpha$, 
$\delta_1,..., \delta_{k-t} \in [n]-\alpha$  
$$ D_{\alpha_1,..., \alpha_t, \delta_1,...., \delta_{k-t}} (T) = 
D_{\alpha_1,..., \alpha_t, \delta_1,...., \delta_{k-t}} $$ 
We prove this on induction on $t$.  
The case $t=0$ is obvious and the case $t=1$ is similar
 to the analogous case 
in the proof of Theorem  \ref{anyk}.
%\medskip
%$ t=1$:  $$ D_{\alpha_j,\delta_1,...., \delta_{s}} (T) \stackrel{de\!f.\; o\!f \; T}{=} 
%a_{\alpha_j} +D_{\underline{\alpha}, \delta_1,...., \delta_{s}} (R) =
%a_{\alpha_j} +D_{\underline{\alpha}, \delta_1,...., \delta_{s}} 
%\stackrel{de\!f.\; of \;D_{\underline{\alpha} , \delta_1,...., \delta_{s}}}{=}$$ 
%$$ = a_{\alpha_j} +D_{\alpha_j, \delta_1,...., \delta_{s}} -a_{\alpha_j} =
%D_{\alpha_j, \delta_1,...., \delta_{s}} $$ 
%\medskip
As to the induction step, suppose first $k-t \geq 1$ 
$$ D_{\alpha_1,..., \alpha_t,\delta_1,...., \delta_{k-t}} (T) 
\stackrel{de\!f.\; o\!f \; T}{=} 
D_{\alpha_2,..., \alpha_t,\delta_1,...., \delta_{k-t}, \delta_{k-t}} (T)
+a_{\alpha_1} -w (\delta_{k-t}, N(\delta_{k-t}, T))=$$  
$$ = D_{\alpha_2,..., \alpha_t,\delta_1,...., \delta_{k-t}, \delta_{k-t}} 
+a_{\alpha_1} -a_{\beta_j} + a_{\beta_j} -
w (\delta_{k-t}, N(\delta_{k-t}, T))=$$  
$$ = D_{\alpha_2,..., \alpha_t,\delta_1,...., \delta_{k-t}, \delta_{k-t}} 
+a_{\alpha_1} -a_{\beta_j} + a_{\beta_j} 
- D_{\beta_i , D,B} (R) + D_{\delta , D,B} (R)
  $$ 
for any $B \subset \beta $ and $ D \ni \delta$. If we take $ A= (\alpha_2,
..., \alpha_t)$, $D= (\delta_1, ...., \delta_{k-t})$ in condition 3, 
we get that the number above  is equal to  
$D_{\alpha_1,..., \alpha_t,\delta_1,...., \delta_{k-t}} $.

Suppose now that $ k-t=0$. We have to prove that 
$ D_{\alpha_1,..., \alpha_{k}} (T) =
 D_{\alpha_1,..., \alpha_{k}}$. We can write 
it as $$ D_{\alpha_1^{s_1},..., \alpha_r^{s_r}} (T) =
 D_{\alpha_1^{s_1},..., \alpha_r^{s_r}} $$ with $ \alpha_1, ...., \alpha_r$ 
distinct and $ s_1+ .... +s_r=k$. We can prove it by induction on 
$ s_3+....+s_r$. 

If $ s_3+....+s_r=0$, the statement is obviuos
 by the definitions of the $ a_{\alpha_i}$ and $T$:
$$ D_{\alpha_1^{s_1},\alpha_2^{s_2}} (T) =
  s_1 a_{\alpha_1} + s_2 a_{\alpha_2} = 
D_{\alpha_1^{s_1}, \alpha_2^{s_2}}$$
If $ s_3+....+s_r>0$, we can suppose for instance that $ s_r >0$:
 $$ D_{\alpha_1^{s_1},..., \alpha_r^{s_r}} (T) =
  D_{\alpha_1^{s_1+1},..., \alpha_r^{s_r-1}} (T) +a_{\alpha_r} -
a_{\alpha_1}= D_{\alpha_1^{s_1+1},..., \alpha_r^{s_r-1}}  +
D_{\alpha_r, X} -
D_{\alpha_1,X}$$ for any $X$.
By taking $X= (\alpha_1^{s_1},..., \alpha_r^{s_r-1})$  
we get $ D_{\alpha_1^{s_1},..., \alpha_r^{s_r}} $.

\hfill \framebox(7,7)

\bigskip

%{\bf Open problem.}
%It is natural to ask if we can generalize the results of this paper to 
%graphs.

\bigskip

\end{document}